\documentclass{tran-l}

\newtheorem{theorem}{Theorem}[section]
\newtheorem{lemma}[theorem]{Lemma}

\theoremstyle{definition}

\newtheorem{example}[theorem]{Example}

\newtheorem{proposition}[theorem]{Proposition}
\newtheorem{corollary}[theorem]{Corollary}
\newtheorem{conjecture}[theorem]{Conjecture}

\theoremstyle{remark}
\newtheorem{remark}[theorem]{Remark}

\def\CC{\mathbb{C}}

\def\ZZ{\mathbb{Z}}
\def\CP{\mathbb{CP}}

\def\cal{\mathcal}  
\newcommand{\A}{{\cal A}}
\newcommand{\B}{{\cal B}}

\numberwithin{equation}{section}

\begin{document}

\title[Derivations of sets of multi-points]{Derivations of an effective divisor on the complex projective
line}


\author{Max Wakefield}
\address{Department of Mathematics, University of Oregon, Eugene, OR 97403}
\curraddr{}
\email{mwakefie@math.uoregon.edu}

\author{Sergey Yuzvinsky}
\address{Department of Mathematics, University of Oregon, Eugene, OR 97403}
\curraddr{}
\email{yuz@math.uoregon.edu}
\thanks{Research at MSRI is supported in part by NSF grant DMS-9810361}


\date{7-25-2005}

\dedicatory{}

\begin{abstract}
In this paper we consider an effective divisor on the complex projective line and associate with it the module D consisting of all
the derivations $\theta$ such that $\theta(I_i)\subset
I_i^{m_i}$ for every $i$, where $I_i$ is the ideal of $p_i$.
The module D is graded and free of rank 2;
the degrees of its homogeneous basis, called
the exponents, form an important invariant of the divisor. We
prove that under certain conditions on $(m_i)$ the exponents do not depend
on $\{p_i\}$. Our main result asserts that if these conditions do not hold
for $(m_i)$ then there exists a general position of $n$ points
for which the exponents do not change. We give an explicit formula for them.
We also exhibit some examples of degeneration of the exponents, in particular
those where the degeneration is defined by vanishing of certain Schur
functions.  As application and motivation, we show that our results imply
Terao's conjecture (about the combinatorial nature of the freeness of hyperplane
arrangements) for certain new classes of arrangements of lines in the complex projective plane.
\end{abstract}

\maketitle


\section{Introduction}
\bigskip
Arguably, the most intriguing conjecture in the theory of hyperplane
arrangements
is the Terao conjecture about the combinatorial character of the freeness
of arrangements. The first interesting and open case is formed by
3-arrangements, i.e., arrangements
of projective lines in the complex projective plane $\CP^2$.
Recent progress has been made by M.~Yoshinaga (see Section \ref{terao})
who found a new relation between the freeness of an arrangement $\A$ and
its
restriction $\A_H$ to an arbitrary hyperplane $H\in\A$.
More precisely the hyperplanes of $\A_H$ have natural multiplicities, which
allows one to consider the multi-arrangement $\tilde\A_H$. Then if $\A$
has
rank $\ell+1$ and is
free with the exponents $\{e_0=1,e_1,\ldots,e_{\ell}\}$ then $\tilde\A_H$
must be free with the exponents $\{e_1,\ldots,e_{\ell}\}$. For
arrangements
in $\CP^2$ this necessary condition is also sufficient (even if it is checked
for only
one line). This result brings to light multi-arrangements of points
(i.e., effective divisors) on $\CP^1$.

In the rest of the paper we will consider multi-arrangements\hfill\break
$\tilde\A=\{p_1,p_2,\ldots,p_n\}$ in $\CP^1$ with $n\geq 2$ and
the positive integer multiplicity $m_i$ of $p_i$ ($i=1,2,\ldots,n$).
We usually order the points subject to
$m_1\geq m_2\geq m_3\geq\cdots\geq m_n$ and call
$m=(m_1,m_2,\ldots,m_n)$ the multiplicity vector of $\tilde\A$.
Also we put $\tilde n=\sum\limits_{i=1}^nm_i$.
If $m=(1,1,\ldots,1)$ then $\tilde\A=\A$ and is called a simple arrangement.

Choosing an appropriate
coordinate system on $V=\CC^2$ we can view the symmetric algebra
$S=\mathcal{S}(V^*)$ as the polynomial algebra
$\CC[x,y]$ and fix linear forms
$\alpha_i=x-z_iy\in S$ ($z_i\in\CC$)
such that $ker(\alpha_i )=p_i$ for all $i$. Then the
defining polynomial of $\tilde\A$ is
$Q=\prod\limits_{i=1}^n(x-z_iy)^{m_i}$.

The graded $S$-module of derivations
 $$D(\tilde{\mathcal{A}})=\{ \theta \in Der_{\CC}S\
|\  \theta (\alpha_i)\in \alpha_i^{m_i}S\}$$
is known to be always free (as a reflexive module over a ring of
homological
dimension 2) of rank 2 whence it is freely generated by two derivations
whose
degrees $e_1$ and $e_2$ are uniquely determined by $\tilde\A$. We always assume that
$e_1\leq e_2$ and call the pair
$exp(\tilde{\mathcal{A}})=(e_1,e_2)$ the exponents of $\tilde\A$.
As it was noticed by G.~Ziegler \cite{Zi}
 (who first considered $D(\tilde\A)$) the exponents are not in general
determined by $m$, unlike in the case of simple arrangement where $e_1=1$ and
$e_2=n-1$, (see \cite{OT}).

The main result of this paper (Theorem \ref{main}) is that under certain
conditions on the multiplicity vector $m$
there exists a general position set (more precisely,
 a nonempty set, open in the Zariski
topology) of $n$ points on $\CP^1$ such that
$$exp(\tilde{\mathcal{A}})=\left(\left\lfloor
\frac{1}{2}\sum\limits_{i=1}^nm_i\right\rfloor ,\left\lceil
\frac{1}{2}\sum\limits_{i=1}^nm_i\right\rceil \right)$$
for every $n$-tuple of points in this position.
The proof of the theorem is broken in several steps and occupies
Sections \ref{theorem}-\ref{proof}. The main part of the proof is in
Section
\ref{proof} where we prove that a certain determinant
is not identically zero as a polynomial in $z_i$'s.

Also we
study all the cases where $m$ does not satisfy these conditions. We prove
that in these cases the exponents are
determined by $m$ and give explicit formulas for them
(Section \ref{combinatorial}). In Section \ref{degen}, we
consider
examples of degeneration of the exponents, in particular a case where the
degeneration is defined by the vanishing of Schur functions for rectangular
diagrams. Finally in Section \ref{terao} we recall the Terao conjecture
and
Yoshinaga theorem and exhibit several new classes of 3-arrangements for
which
our results imply the conjecture.

This paper was
partially written when both authors were participating in the MSRI program
on
Hyperplane Arrangements and Applications. We are grateful to MSRI for
support.  We also are grateful to Hiro Terao and Masahiko Yoshinaga for
useful discussions.

\section{Multiplicity vectors that determine exponents}
\label{combinatorial}

First, for convenience of the reader we recall Zieler's generalization (\cite{Zi}, p. 351) of Saito's criterion (\cite{OT}, Theorem 4.19) restricted to our case. Let $\theta_1,\theta_2\in D(\tilde{\mathcal{A}})$ and in some coordinate system on $\CP^1$ we have $\theta_i=f_{i1}\partial_x+f_{i2}\partial_y$ for some homogeneous $f_{ij}\in S$ of the same degree where $i,j=1,2$ and $\partial_x$ and $\partial_y$ are derivative with respect to $x$ and $y$. Then $(\theta_1,\theta_2)$ is a basis of $D(\tilde{\mathcal{A}})$ if and only if the determinant $f_{11}f_{22}-f_{12}f_{21}$ is equal to $\tilde{Q}$ multiplied by a nonzero constant.

Now we mention two simple properties of the exponents$(e_1,e_2),
e_1\leq e_2,$ which we use frequently in the rest of the paper.

1. $e_1+e_2=\tilde n$.

This follows from the version of Saito's criterion of the freeness for
multi-arrangements.

2. Suppose $\A_1\subset\A_2$ are two multi-arrangements and
$\tilde n(\A_2)=\tilde n(\A_1)+1$. If $exp(\A_1)=(e_1,e_2)$ and
$exp(\A_2)=
(e_1',e_2')$ then either $e_1'=e_1$ (whence $e_2'=e_2+1$) or $e_2'=e_2$
(whence $e_1'=e_1+1$).

Indeed since $D(\A_2)\subset D(\A_1)$ we have $e_i'\geq e_i$, $i=1,2$.
Now the property follows from property 1.

In the rest of the section we will
consider several cases of multiplicity vectors
$m$ where the exponents $(e_1,e_2)$
 are uniquely determined by $m$ and do not depend on the
position of the points. In most of these cases we will be able to
 exhibit a derivation of minimal degree from
$D(\tilde\A)$ as a function of $m$.

In the rest of the paper we always assume that a coordinate system is
fixed
on $\CP^1$ so that the defining polynomial of $\tilde\A$ is
$$\tilde Q=x^{m_1}y^{m_2}\prod\limits_{i=3}^n(x-z_iy)^{m_i}\eqno(2.1)$$
for some $z_i\in\CC\setminus\{0\}$. Then any derivation from $D(\tilde\A)$
has the form $fx^{m_1}\partial_x-gy^{m_2}\partial_y$ for some $f,g\in
S$
($\partial_x$ and $\partial_y$ are derivatives with
respect to $x$ and $y$).

We can exhibit two derivations in $D(\tilde \A)$ for future use.
Put $\theta_1=\frac{\tilde{Q}}{x^{m_1}}{\partial_y}$. It is easy to see
that $\theta_1$ is a homogeneous derivation having minimal degree
among
all the elements of $D(\tilde\A)$ with $f=0$. Also put
$\theta_2=\frac{\tilde{Q}}{Q}\theta_E$, where $\theta_E=x
{\partial_x}+y{\partial_y}$ is the Euler derivation. Again it is easy to
check that $\theta_2$ is a homogeneous derivation having minimal
degree
among all elements of $D(\tilde\A)$ proportional to $\theta_E$.

{\bf Case 2.1}. Let $m_1 \geq \sum
\limits_{i=2}^nm_i$.
 Then $\theta_1$ is a derivation of minimal degree in $D(\tilde\A)$.
In particular $(e_1,e_2)=(\sum\limits_{i=2}^nm_i,m_1)$.

\proof Since any homogeneous derivation
$\theta$ from $D(\tilde\A)$ with $f\not=0$ satisfies
$deg(\theta ) \geq m_1 \geq \sum \limits_{i=2}^n
m_i=deg(\theta_1 )$, $\theta_1$ has the smallest degree
in $D(\tilde{\mathcal{A}})$ and the statement follows.
\qed

{\bf Case 2.2.} If $\tilde n\leq 2n-2$
then $\theta_2$
is a derivation of minimal degree in $D(\tilde{\mathcal{A}})$.
In particular
$$(e_1,e_2)=(\tilde n-n+1,n-1).$$

\proof
 Notice $deg(\theta_2)= \tilde n-n+1\leq 2n-2-n+1=n-1$. Suppose that there
is
$\theta\in D(\tilde{\mathcal{A}})$ such that $deg(\theta)
<deg(\theta_2)\leq n-1$.
Since $D(\tilde{\mathcal{A}})\subset D({\mathcal{A}})$ and the latter
module
is generated by $\theta_E$ together with a derivation of degree $n-1$ we have that
$\theta=h\theta_E$ for a polynomial $h$. Applying $\theta$
to $\alpha_i$ we see
that $h$ is divisible by $\frac{\tilde{Q}}{Q}$ which contradicts
the
condition on degree of $\theta$.

Hence the degree of $\theta_2$ is minimal and the result follows.
\qed

\begin{remark}
\label{odd}
(i) Cases 2.1 and 2.2 meet
 at the case $m=(n-1,1,\ldots,1)$.
For this $m$ the derivations found in the cases give a minimal generating
system of degrees $(n-1,n-1)$.

Indeed assume $\tilde n\leq 2n-2$
and $m_1\geq\sum\limits_{i=2}^nm_i$. Since every $m_i\geq 1$
 we have $\sum\limits_{i=2}^{n}m_i$ and $m_1$ both are not less than
$n-1$ whence
 $\tilde n=2n-2$, $m_1=n-1$ and $m_i=1$ for every $i>1$.

(ii) The previous remark allows us to handle the case where $\tilde n=2n-1$.
Since $n\geq 2$ there exists an $m_i>1$. Decreasing this $m_i$ by 1 we obtain
a
new multi-arrangement with exponents $(n-1,n-1)$. By property 2, the
exponents
of the initial arrangement are $(n-1,n)$. In particular they are
determined
by $m$.

\end{remark}
Together with Case 2.2, these comments imply that  if the average of
all $m_i$ is less than 2
then the exponents are determined by $m$. If the average equals 2, this is
false.
We however can give one example of $m$ for which it is true.

\begin{example} If $m_i=2$ for all $1\leq i\leq n$ then the
exponents are
$(n,n)$. Indeed we can exhibit two derivations of degree $n$ satisfying
the Saito criterion.

Let an arrangement $\A$ be defined by
$Q=xy\prod\limits_{i=3}^n(x-z_iy)$ for some $n-2$ numbers $z_i\in\CC$
and $\tilde Q=Q^2$. If $n=2$ we can take $\theta_1=x^2\partial_x$ and
$\theta_2=y^2\partial_y$ and the result obviously follows.

 Suppose $n>2$ and for every $i,\ 3\leq i\leq n,$ put
$Q_i=\frac{Q}{x-z_iy}$, $h_1=\sum\limits_{i=3}^n Q_i$, and
$h_2=\sum\limits_{i=3}^n z_iQ_i$. Then put
$$\xi_1=h_1\theta_E+\frac{y}{x}Q\partial_y$$
and
$$\xi_2=h_2\theta_E-\frac{x}{y}Q\partial_x.$$
One can check straightforwardly that $\xi_i\in D(\tilde \A)$ ($i=1,2$)
and the matrix of the coordinates of
these derivations has determinant equal to
$(n-2)Q^2$.

Let us remark that just a slightly different basis for this case was previously constructed in \cite{ST}, Proposition 5.4.
\end{example}

If $n=2$ or $3$ then of course the exponents are determined by $m$ since
there is only one position of 2 or 3 points on $\CP^1$ up to a projective
isomorphism. For $n=4$ even for the average of $m_i$ equal 2
the situation is different. It follows from the cases above that the vectors
$(5,1,1,1)$, $(4,2,1,1)$, and $(2,2,2,2)$ determine exponents uniquely
(they are $(3,5)$ for the first vector and $(4,4)$ for two others).
On the other hand, for the vectors $(3,3,1,1)$ and $(3,2,2,1)$ the exponents depend
on
the position of 4 points on $\CP^1$ (see section \ref{degen}).

\section{The main theorem}
\label{theorem}

In this section, we are concerned with an arbitrary
 multiplicity vector $m$ not satisfying the strict inequalities
studied in the previous section.
Our main result asserts that $m$ uniquely determines the exponents
$(e_1,e_2)$ for $n$ points in general position on $\CP^1$.
(In this paper we understand by a general position set of $n$ points
in $\CP^1$, a nonempty set in the configuration space
$(\CP^1)^n$ open in the Zariski topology.)
Roughly speaking these general position exponents are as close to being
equal as possible.
 More precisely we are going to prove the following.

\begin{theorem}
\label{main}
Suppose $m$ has the following properties: $m_1\leq\sum\limits_{i=2}^nm_i$
 and $\tilde n\geq 2n-2$. Then there exists a general position
of $n$ points such that
$$exp(\tilde{\mathcal{A}})=\left(\left\lfloor
\frac{1}{2}\tilde n\right\rfloor ,\left\lceil
\frac{1}{2}\tilde n\right\rceil \right)$$
for every multi-arrangement $\tilde\A$ of points in this position
having $m$ as the multiplicity vector.
\end{theorem}

The proof will be broken into several parts according to the following plan.
First we focus on the case when $\tilde n$ is even. Then for
$\tilde n$ odd the result can be deduced similarly to Remark \ref{odd} (i)
(cf. section \ref{proof}).
For $\tilde n$ even the condition guaranteeing the existence of a derivation in
$D(\tilde\A)$ of degree less than $\tilde n/2$ can be naturally expressed
as
the vanishing of the determinant of a square matrix $M$ of size
$\sum\limits_{i=3}^nm_i$ whose entries are monomials in $z_i$,
$i=3,4,\ldots,n$.
 Thus the determinant of $M$ is a polynomial, $d$, in these indeterminates and
to prove the theorem it suffices to show that $d$ is not
identically 0. For this, in turn, it suffices to show that the leading
monomial
of $d$  (in some linear ordering of the monomials) has a nonzero
coefficient.

Now the proof branches out. Since the set of rows of $M$ is by
construction
partitioned into blocks corresponding to the points $p_3,p_4,\ldots,p_n$,
 it is natural
to use the Laplace formula for $d$. In most cases the coefficient of the
leading term is the product of minors of $M$, one from each block of
rows, and it is not hard to show that none of these minors vanishes.
However, there are
cases where the leading term is not a single product but a sum of
several products of
minors. This case is harder and requires a deeper analysis of the
minors.

\section{The matrix $M$}
\label{matrix}

In this section, we assume that $\tilde n$ is even
and put $e=\frac{1}{2}\tilde n-1$. We assume also that $m_i\geq 2$ for
every
$i$ and
$m_1<\sum\limits_{i=2}^nm_i$. Equivalently, $e-m_1\geq 0$ whence
also
$e-m_2\geq 0$. Any derivation $\theta\in
D(\tilde{\mathcal{A}})$ of degree $e$
is of the form $\theta
=x^{m_1}f(x,y)\partial_x-y^{m_2}g(x,y)\partial_y$ where $f(x,y),g(x,y)\in
S$, $deg(f(x,y))=e-m_1$ and $deg(g(x,y))=e-m_2$. Write
$$f(x,y)=\sum\limits_{j=0}^{e-m_1}f_jx^{e-m_1-j}y^j$$
and
$$g(x,y)=\sum\limits_{j=0}^{e-m_2}g_jx^{e-m_2-j}y^j$$ where
$f_j,g_j\in \CC$. Besides since $\theta\in
D(\tilde{\mathcal{A}})$, for all $i\geq 3$ we have
$\theta(x-z_iy)\in (x-z_iy)^{m_i}S$.
Sacrificing homogeneity, we put $y=1$ and obtain the equivalent
condition
$$\partial_x^k(\theta(x-z_iy)|_{y=1})|_{x=z_i}=0 \eqno(*)$$
 for all $0\leq k\leq m_i-1$
where $\partial_x^k$ is the $k^{th}$ derivative with respect to $x$.
Treating the coefficients $f_i$ and $g_j$ as unknowns we consider (*)
as a system of equations for these unknowns. Explicitly the equation
corresponding to the point $p_i\ (3\leq i\leq n)$ and some $k$
($0\leq k\leq m_i-1$) is
$$\sum\limits_{j=0}^{e-m_1}f_j\frac{(e-j)!}{(e-j-k)!}
z_i^{e-j-k}+z_i\sum\limits_{j=0}^{e-m_2}
g_j\frac{(e-m_2-j)!}{(e-m_2-j-k)!}z_i^{e-m_2-j-k}=0$$
where we agree that $z_i^{\ell}=0$ for $\ell<0$.

The matrix $M$ is the matrix of coefficients of this
system of equations which are polynomials in $z_i$, $i=3,4,\ldots,n$.
A simple computation shows that the matrix is square - the number of rows
as
well as columns is $\sum\limits_3^nm_i$. Moreover the set of rows of $M$
is partitioned into $n-2$ blocks $L_3,L_4,\ldots,L_n$ where the block
$L_i$
does not contain $z_j$ with $j\not= i$.
The set of columns of $M$ is partitioned into
two blocks:  the $f$-block and the $g$-block
consisting of the coefficients at $f_i$ and $g_j$
respectively. The sizes of the blocks are $e-m_1+1$ and $e-m_2+1$
respectively. Notice that the $f$-block of the block $L_i$ is a Wronskian
matrix of the functions $x^e,x^{e-1},\ldots,x^{m_1}$ for their derivatives
up
to $m_i-1$ evaluated at $z_i$, i.e., an entry of the block has the form
$\partial_x^{k}(x^{e-\ell})|_{x=z_i}$,
$k=0,1,\ldots,m_i-1,\ell=0,1,\ldots,
e-m_1$. The $g$-block of the block $L_i$ is a Wronskian matrix of the
functions
$x^{e-m_2},x^{e-m_2-1},\ldots,1$ for the same range of derivatives as for
the $f$-block, multiplied by $x$ and then evaluated at $z_i$. So its
typical
entry has the form $(x\partial_x^k(x^{e-m_2-\ell}))|_{x=z_i}$,
$k=0,1,\ldots,m_i-1,\ell=0,1,\ldots,e-m_2$.

By construction of $M$ the vanishing of its determinant at some $(n-2)$
-tuple $(z_3,\dots,z_n)$ of pairwise different nonzero
complex numbers is equivalent to the existence of a derivation of degree
$e$
in $D(\tilde\A)$  where $\tilde\A$ is the multi-arrangement defined by
the formula (2.1).
 We can view $z_i$ as
indeterminates and $\det(M)$ as a polynomial in these $n-2$ indeterminates.
In the space
$\CC^{n-2}$ of  all $n-2$-tuples of points the complement
to the zero locus of $d$ is an
open set in the Zariski topology. Thus Theorem \ref{main} under the conditions on $m$
assumed
in this section would follow from the following result.

\begin{theorem}
\label{det}
The polynomial $d$ is not identically 0.
\end{theorem}
The rest of the section is devoted to the proof of Theorem \ref{det} which
we break into several parts.

\begin{lemma}
\label{vander}
The Wronskian of the power functions $x^{\lambda_1},\ldots,x^{\lambda_k}$
where $(\lambda_i)$ is a strictly decreasing $k$-tuple of nonnegative
integers
is a monomial in
$x$ of degree $\sum\limits_{j=1}^k\lambda_j-\sum\limits_{r=1}^{k-1}r$
with
the coefficient
$$(-1)^{\lfloor\frac{k}{2}\rfloor}\prod\limits_{1\leq i<j\leq
k}(\lambda_i-\lambda_j).$$
\end{lemma}

\proof  The first part of the statement
follows immediately from the fact that every summand in the general
formula
for the determinant is the product of exactly one entry from each row and
each column.

The second part of the statement can be
 proved by setting $x=1$ and transforming the obtained matrix
to a Vandermonde matrix by applying row operations and
induction
on number of rows. The sign
$(-1)^{\lfloor\frac{k}{2}\rfloor}$ comes from the unusual ordering of the
rows of the Vandermonde determinant.
\qed

In order to apply Lemma \ref{vander} we represent $d$ via (repeated) Laplace's development
corresponding to the fixed partition
$(L_3,\ldots,L_n)$ of the set of rows of $M$. For that we need to consider ordered partitions
$\beta=(B_3,\ldots,B_n)$ of columns such that $|B_i|=|L_i|=m_i$ for every $i$. Each $\beta$ defines a
permutation $P(\beta)$ of all the columns of $M$. This gives a bijection of the set of these partitions
with the set of permutations of the columns of $M$ such that the orderings of columns
inside subsequent intervals, $B_i$, of length $m_i$ are induced from $M$.
Then we have
$$d=\sum\limits_{\beta}{\rm sign}(P(\beta))\prod\limits_{i=3}^n d(L_i,B_i)$$
where the summation is over all the ordered partitions $\beta$ as above
and $d(L_i,B_i)$ is the
$m_i\times m_i$-minor of
$M$ with rows from $L_i$ and columns from $B_i$.

\begin{remark}
The usual formula for signs of summands in Laplace's development differs from ours.
Our formula can be easily checked by considering the signs of products of diagonal elements of each
minor $d(L_i,B_i)$$($see \cite{Mu}$)$.
\end{remark}

The degrees of entries in each column are determined by the degree $r$ of
the top element in the column; for brevity we say
that {\it the column has degree $r$}.

In order to prove that the polynomial $d$ does not vanish identically it
suffices to find a monomial of it with a nonzero coefficient. We will use
for this purpose the leading term of $d$ in the lexicographic order
generated by the reverse ordering of
the indeterminates  $z_n<z_{n-1}<\cdots<z_3$.
The proof of Theorem \ref{det} is immediate when the leading term is the
product of (uniquely determined) minors, one from each block $L_i$.
The simplest case when this
happens is the case of `no overlapping'.

\begin{proposition}
\label{nooverlap}
Suppose that the sets of  degrees of (the top rows of)
columns in
the $f$-block and those in the $g$-block are disjoint (equivalently
$m_1+m_2> \frac{\tilde n}{2}$). Then Theorem \ref{det} holds.
\end{proposition}
\proof  Indeed in this case there are no two columns of same degrees.
There is the unique partition $\beta$ of the columns such that
the corresponding product of minors is the leading term of $d$;
it is the partition corresponding to the identity permutation of the columns.
The respective minors
do not vanish by Lemma \ref{vander}.    \qed

In the rest of the proof we focus on the case where there is some
overlapping.  This means that several columns on the right flank of
the $f$-block have the same
degrees as several columns on the left flank of the $g$-block.
The number of columns of this
kind in the $f$-block is the same as this number for the $g$-block and
equals $s=\frac{1}{2}\tilde n+1-(m_1+m_2)$. We denote by $O$ the set of
 these columns and denote its columns by
$a_1,\ldots,a_s$ and $b_1,\ldots,b_s$ meaning that the degrees of $a_i$
and
$b_i$ are the same ($i=1,2,\ldots,s$) and $deg( a_i)=deg( a_{i+1})+1$
($i=1,2,
\ldots,s-1$). Notice that $deg(a_1)=deg( b_1)=m_1+s-1=e-m_2+1$ and
the number of columns with higher degrees is $e-deg( a_1)=m_2-1$.
The only minimal
linearly dependent, over $S$, sets of columns in any block $L_j$ are $\{a_{i+1},b_{i}\}$
$(i=1,2,\ldots,s-1)$.

Even when overlapping occurs there still may be
a unique choice of minors from the $L_i$'s whose product
equals the leading term of $d$. However, the simple partition
that worked in the non-overlapping case may not work now because some minors may vanish.

To analyze the situation deeper we need more notation.
Consider a partition $\beta=(B_3,\ldots,B_n)$ of columns such that the corresponding product of minors
equals the leading monomial of $d$ and the respective permutation $P(\beta)$.
Then $P(\beta)$ is the identity permutation on the columns outside of $O$.

To express this in terms of partitions denote by $B_{i_1}$ and $B_{i_2}$ the first and the last blocks
intersecting with $O$.
Then the blocks $B_3,\ldots,B_{i_1-1}$ and $B_{i_2+1},\ldots,B_n$
  are determined uniquely
 as well as $r_1$ columns in
$B_{i_1}$ and $r_2$ columns in $B_{i_2}$ not in $O$
(for some $0\leq r_j < m_{i_j}$, $j=1,2$).
For convenience, put $m_i'=m_i$ if $i_1<i<i_2$ and $m_{i_j}'=m_{i_j}-r_j$
for $j=1,2$. In particular $\sum\limits_{i=i_1}^{i_2}m'_i=2s$.

Now we can describe a (perhaps not unique) partition of $O$
recursively.

\begin{proposition}
\label{minors}
Let $i_1\leq i\leq i_2$. Suppose $j$ is the minimal index such that
the column $a_j$ is not chosen for all blocks $B_{i'}$ with $i'<i$
and $k$ is similar index for $b_k$.

If $j\leq k$ then the block $B_i$
has the following columns from $O$:

$$\begin{array}{ll}
(1)\  a_j,a_{j+1},\dots,a_{k-1},a_k,b_k,b_{k+1},\ldots,b_{m'_i+j-2}
 & if\ m'_i\geq k-j+2\\
(2)\  a_j,a_{j+1},\ldots,a_{m'_i+j-1} & if\ m'_i\leq k-j\\
(3)\  a_j,a_{j+1},\ldots,a_{k-1},b_k & if\ m'_i=k-j+1\\
(4)\  a_j,a_{j+1},\ldots,a_{k-1},a_k & if\ m'_i=k-j+1\\
                               & {\rm and\ either }\ i=i_2\
                               {\rm or}\ m'_{i+1}\leq 2.\\
\end{array} $$

If $j>k$ then $m'_i=2$, $j=k+1$ and the columns in $B_i\cap O$ are
$$b_k,b_j.$$

\end{proposition}
\proof
We use induction on $i$. It is obvious for $i=i_1$ and follows from the
inductive hypothesis for $i>i_1$ that the indexes of columns $a_t$ and $b_t$
available for $B_i$ form the intervals $j\leq t\leq s$ and
$k\leq t\leq s$ respectively. Choosing the columns of maximal available
degrees and avoiding linear dependent pairs we immediately obtain the
first
three cases for $j\leq k$.

The fourth case is more subtle. If $m'_i=k-j+1$
then the columns
from $a_j$ to $a_k$ have the maximal degrees and are independent
whence seem to be appropriate for $B_i$.
Assume however that $i\not=i_2$ and $m'_{i+1}\geq 3$. Then with the above
choice of columns for $B_i$ the block $B_{i+1}$
has columns
$b_k,b_{k+1},b_{k+2},\ldots$. However the choice of columns for
$B_i$ from
the third case gives the same degree minor in $L_i$ and allows the higher
degree choice $a_k,a_{k+1},b_{k+1},\ldots$ for a minor in $L_{i+1}$.
Thus, the choice (4) in this case does not work.
In the case where $m'_{i+1}=1$ or $2$ this obstruction disappears and both
choices (3) and (4) for $B_i$ give nonzero minors of the same degree.

Finally consider the case $j>k$. This condition means that the column
$a_k$ was used in the previous blocks and $b_k$ was not. Applying the
inductive hypothesis we can see that it is possible if and only if the
choice (4) occurs in $L_{i-1}$ whence $m'_i=2=k-j+1$.
The choice of columns for $B_i$ is now clear.
\qed

\begin{remark}
Proposition \ref{minors} can be used to show that even if $O\not=\emptyset$
there may be just one product of minors equal to the leading term of $d$.
Moreover, it is possible in general to compute the index $r$ ($i_1\leq r\leq
i_2$) of the block, where the non-uniqueness
starts, directly from the vector $m$.
Since we won't use this in the paper we omit this computation.
\end{remark}
We will need however some properties of $r$ which we collect in the
following corollary that
follows immediately from Proposition \ref{minors}.

\begin{corollary}
\label{notunique}
Suppose that for a multiplicity vector $m$ there are at least two different
partitions of $O$ whose corresponding products of minors equal the leading monomial of $d$.
Then

(1) The union of $\bigcup\limits_{r'<r}B_{r'}$ and first $m_r-1$ columns
of $B_r$ consists of all the columns of degree larger than $deg( a_t)$
for some $t\leq s$ satisfying the equality
$$2(t-1)+m_2=\sum\limits_{i=3}^rm_i.$$

(2) Either $t=s$ or $m_i=2$ for all $i>r$.

(3) The set of columns
$O'=\{a_t,a_{t+1},\ldots,a_{s},b_t,b_{t+1},\ldots,b_{s}\}$
is partitioned by the blocks from $B_{r}$ to $B_{i_2}$.
This partition has the following properties:

(a) the block $B_r$ has only one column from $O'$, either $a_t$ or $b_t$
(which implies that $B_{i_2}$ has either $a_s$ or $b_s$ and no other columns
from $O$);

(b) if the second column of a block
is $a_j$ then the first
column of  the next block is $b_j$ and vice versa;

(c) if the first column of a block is $a_j$ then its second
column is either $a_{j+1}$ or $b_{j+1}$;

(d) if the first column of a block is $b_j$ then its second column
is $b_{j+1}$.

\end{corollary}

The ordered partitions of $O'$ having the properties of the
previous Corollary are in one-to-one
correspondence with some permutations of $O'$ that we call admissible.

\begin{example}
\label{3}
Using the notation of the previous Corollary, if $k=1$ and $s=2$ then
there exists the following 3 admissible partitions of $O=O'$
$$a_1|b_1b_2|a_2,\ b_1|a_1a_2|b_2,\ b_1|a_1b_2|a_2.$$
\end{example}

It is easy to explicitly compute minors involved in Corollary
\ref{notunique}.
Under the conditions of Proposition \ref{minors}
consider any product of minors giving maximal degree monomial of $d$.
Recall that the minor $\mu_i=d(L_i,B_i)$
 is a monomial in $z_i$ with some coefficient
$c_i\in\ZZ, \ i=3,\ldots,n$.

\begin{lemma}
\label{1-2}
The coefficients, $c_i$, are as follows.

If $i_2-r>1$ then

(i) $c_i=-1$ for $r<i\leq i_2$ if
both columns of $\mu_i$ are from the
$f$-block or both columns are from $g$-block and
$c_i=-2$ if the columns are from different
blocks.

If $i_2=r+1$ then

(ii)  $c_{i_2}=(-1)^{\lfloor\frac{m_{i_2}}{2}\rfloor}\pi(m_{i_2}-1)!$
if the first column of $\mu_{i_2}$ is $b_s$ and
$c_{i_2}=(-1)^{\lfloor\frac{m_{i_2}}{2}\rfloor}\pi m_{i_2}!$ if the first column is $a_s$
where $\pi=\prod\limits_{j=1}^{m_{i_2}-2}j!$ for
$m_{i_2}>2$ and $\pi=1$ otherwise.

For arbitrary $i_2-r\geq 1$

(iii) $c_{r}=(-1)^{\lfloor\frac{m_r}{2}\rfloor}\pi'(m_{r}-1)!$
if the last column of $\mu_{r}$ is $a_k$ and
$c_{r}=(-1)^{\lfloor\frac{m_r}{2}\rfloor}\pi' m_{r}!$ if the last column is $b_k$
where $\pi'=\prod\limits_{j=1}^{m_{r}-2}j!$
for $m_r>2$ and $\pi'=1$ otherwise.

\end{lemma}

This lemma follows immediately from Corollary \ref{notunique} and
Lemma \ref{vander}.

Now we are ready to prove that even when the leading monomial of the
determinant $d$ is not just a product of minors but the sum of several of those,
its coefficient does not
vanish. We prove this by computing that coefficient.

For each admissible permutation $P$ of $O'=(a_k,\ldots,a_s,b_k,\ldots,b_s)$
(i.e., a permutation that defines a partition satisfying Corollary
 \ref{notunique})
consider the product $\prod(P)=\prod c_i$ where $r\leq i\leq i_2$ and
 $c_i$ is defined for the block $B_i$ in Lemma \ref{1-2}.
 Notice that the product $\pi'(m_{r}-1)!$ can be factored out
of $\prod(P)$  for each $P$ along with the product
$\prod\limits_{i=3}^{r-1}c_i\prod\limits_{i=i_2+1}^nc_i$.
More precisely
for an arbitrary admissible permutation $P$ put
$$
\prod(P)=((-1)^{\lfloor{\frac{m_r}{2}}\rfloor + \lfloor\frac{m_{i_2}}{2}\rfloor}
\pi'(m_r-1)\prod\limits_{i=3}^{r-1}c_i\prod\limits_{i=i_2+1}^nc_i)C(P)$$
for some $C(P)\in \CC$.

\begin{proposition}
\label{coefficient}
If $i_2=r+1$ then
\begin{equation}\label{4.1}
\sum\limits_P {\rm sign}(P) C(P)=\pi(m_{i_2}-1)!(1-m_rm_{i_2});
\end{equation}

if $i_2>r+1$ then
\begin{equation}\label{4.2}
\sum\limits_P {\rm sign}(P) C(P)= (-1)^{\lceil{\frac{u+2}{2}}\rceil}((m_r-1)u+1)
\end{equation}

 where the sum is taken over all
the partitions $P$ of $O'$ satisfying
Corollary \ref{notunique} and $u=i_2-r+1$.
\end{proposition}

\proof
(i) Suppose $i_2=r+1$.
 Lemma \ref{1-2} (ii) and (iii) implies that there are only two
admissible permutations. The formula (4.1) follows also from that lemma.

(ii) Suppose that $i_2>r+1$, i.e., $u>2$. Corollary \ref{notunique}
implies that $m_{i_2}=2$ whence the left hand side of (4.2) depends only on
$m_r$ and $u$. We denote it $\sigma(m_r,u)$ and put
$\sigma(u)=\sigma(2,u).$

Let us consider first the most important case where $m_r=2$.
The key observation here is the following recursive formula:
$$\sigma(u)= \sigma(u-2)+(-1)^{u}2\sigma(u-1)$$
for every $u\geq 4$. Consider first all permutations $P$ with $a_k$ as the
first element. Then the beginning of $P$
is $a_k|b_{k}b_{k+1}|a_{k+1}$. That shows that the
the $C(P)=C(P')$
where $P'$ is the restriction of $P$ to the rest of $O'$. Observe
that  the size of the set is $u-2$
 and ${\rm sign}(P')={\rm sign}(P)$. This explains the
first summand of the formula.

Then consider all permutations $P$ with $b_k$ as the first element.
This defines the beginning of $P$ as $b_k|a_k$. Thus $C(P)=(-2)C(P')$
where $P'$ is again the restriction of $P$ to the rest of $O'$, this time
of size $u-1$.
Besides ${\rm sign}(P')=(-1)^{u-1}{\rm sign}(P)$ as it is easy to see
by counting transpositions. That explains the second summand.

Augment the recursive formula
 by the initial conditions $\sigma(2)=3$ and $\sigma(3)=-4$
(for $u=3$ there are 3 admissible orderings giving the summands equal
to 2,2,-8, cf. Example \ref{3}).
Using the recursive formula and the initial conditions,
formula (4.2) for $m_r=2$
can be proved by induction on $u$.

If $m_r>2$ then the recursive formula is a little different;
$$\sigma(m_r,u)=\sigma(u-2)+(-1)^{u}m_r\sigma(u-1)$$
although the proof is the same.
Substituting (4.2) for $m_r=2$ in this formula we complete the proof.
 \qed

Proposition \ref{coefficient} completes the proof of Theorem \ref{det}.

\section{Proof of Theorem \ref{main}}
\label{proof}

In this section we finish the proof of Theorem \ref{main}.

Recall that we consider multi-arrangements $\tilde{\mathcal{A}}$
with the following conditions on the multiplicity vector
$(m_1,m_2,\ldots,m_n)$  such that $m_1\geq m_2\geq\cdots\geq m_n$:

(i)$m_1\leq\sum\limits_{i=2}^nm_i$;

(ii)$\tilde n=\sum\limits_{i=1}^nm_i\geq 2n-2$.

If $\tilde n$ is even then condition (i) allows us to write
a square matrix $M$ whose determinant $d=\det M$ is a polynomial in
the coordinates of the points of $\A$. The equation $d=0$ describes the locus
in
$(\CP^1)^{n-2}$ of all arrangements $\A$ such that $e_1(\tilde\A)<\frac{n}{2}$.
For the case $m_n>1$ (i.e., $m_i>1$ for
every $i$), we proved in the previous section that $d$ is not identically
0
which implies the statement of the theorem for this case.

Now we consider the general case where multiplicity vector is arbitrary
(satisfying (i) and (ii)). We apply
induction on the number, $\alpha$, of multiplicities equal 1 using as the
base results of the previous section for $\alpha=0$. Suppose $\alpha>0$
(in particular $m_n=1$)
and consider $m'=(m_1,\ldots,m_{n-1})$.
 If $\tilde n$ is odd then the result
follows immediately from the inductive hypothesis for $m'$, the
monotonicity of $e_1$, and the fact that the pre-image of a general
position
set under a coordinate projection is in general position.
Thus we can assume that $\tilde n$ is even.

Assume that
the result does not hold for $m$. Since $\tilde n$ is even we can again
bring the matrix $M$ into consideration. The assumption implies that $d$ is
the
0 polynomial, i.e., every arrangement with the multiplicity vector $m$ has
\begin{equation}\label{5.1}
e_1\leq \frac{\tilde n}{2}-1.
\end{equation}
In particular if we fix a multi-arrangement $\tilde\A'$ of $n-1$ points in general
position with the multiplicity vector $m'$
then adjoining any extra point $q$
( with multiplicity 1) we obtain a multi-arrangement satisfying (5.1).
Applying the inductive hypothesis to $m'$ we obtain that
for $\tilde\A'$
we have
$\exp(\tilde\A')=(\frac{\tilde n}{2}-1,\frac{\tilde n}{2})$.
 If $\theta$ is
a (unique up to a nonzero multiplicative constant)
derivation from $D(\tilde\A')$ of the degree  $\frac{\tilde n}{2}-1$ then for
any point, $q$, $\theta$ will preserve $q$. Thus, $\theta$
is in $D(\B)$ where $\B$ is an arbitrary (simple) arrangement.
 Choose $\B$ of a
cardinality greater than  $\frac{\tilde n}{2}+1$. Then
$\exp(\B)=(1,e)$ where the exponent 1 corresponds to the
Euler derivation $\theta_E$ and $e>\frac{\tilde n}{2}$.
Hence, $\theta$ is proportional to
$\theta_E$,
i.e., $\theta=P\theta_E$ for a homogeneous polynomial $P$. Since besides
$\theta\in D(\tilde\A)$ we have $deg( P)\geq \tilde n-n$ (see Section
\ref{combinatorial}) which implies
\begin{equation}\label{5.2}
\tilde n-n+1\leq \frac{\tilde n}{2}-1.
\end{equation}
The inequality (5.2) is equivalent to $\tilde n\leq 2n-4$ which contradicts
the condition (ii) of the theorem.
This contradiction completes the proof.         \qed

\bigskip

\section{Degeneration of the exponents}
\label{degen}

In this section, we consider several cases of degeneration of exponents,
i.e., examples of (classes of) multi-arrangements satisfying the
conditions
of Theorem \ref{main} and having
$e_1<\frac{\tilde n}{2}$. We assume for convenience that $\tilde n$ is
even and as before represent the defining polynomial of a
multi-arrangement as $\tilde Q=x^{m_1}y^{m_2}\prod\limits_{i=3}^n
(x-z_iy)^{m_i}$
($m_i\leq m_j$ for $j<i$).
It is necessary for degeneration that the $(n-2)$-tuple $(z_3,\ldots,z_n)$
annihilates the polynomial $d$.

\begin{remark}
The polynomial $d$ has a factor whose zero locus does not contain
acceptable $n-2$-tuples. Factoring $z_i$ from each row of a block $L_i$
of the matrix $M$,
then subtracting similar rows of different blocks and factoring
differences
$(z_i-z_j)$ we obtain
$$d=\prod\limits_{k=3}^nz_k^{m_k}\prod\limits_{3\leq j<i\leq
n}(z_j-z_i)^{m_i}d_1$$
for some polynomial $d_1$.

 We conjecture that $d_1$ is not divisible by any
$z_i$ or $z_i-z_j$. There could be however $(n-2)$-tuples in the zero
locus
of $d_1$ that do not define an arrangement (with $n$ points).
For instance, if $m=(6,4,3,2,1)$ then $d_1(1,1,1)=0.$
\end{remark}

{\bf Examples.}

1. $m=(3,3,1,1)$. This example was considered first by G.Ziegler
\cite{Zi}. An arbitrary multi-arrangement with the multiplicity vector $m$
is projectively equivalent to $x^3y^3(x-y)(x-zy)$
 ($z\in\CC\setminus\{0,1\}$) and the matrix $M$ is
$$M=\left (\begin{matrix}
    1 & -1\\
    z^3 & -z
\end{matrix} \right ).$$
The polynomial $d_1=z+1$ vanishes at $z=-1$.
Thus the only degenerate case is $x^3y^3(x^2-y^2)$ that have exponents
(3,5). Notice that the degeneration here happens exactly in the case where
the quadruple $(p_1,p_2,p_3,p_4)$ of points is harmonic (i.e., its
cross-ratio equals -1).

2. An interesting generalization of the previous example is given by
the vector $m=(k,k,1,\ldots,1)$ where $k\geq 3$ and $n=k+1$. A respective
multi-arrangement can be written as
$x^ky^k\prod\limits_{i=1}^{k-1}(x-z_iy)$.
As a straightforward computation shows this arrangement has $e_1=k$
(the minimal possible) if and only if the set
$\{z_1,z_2,\ldots,z_{k-1}\}$ equals to the set of all roots of degree
$(k-1)$
of a nonzero (complex) number. Notice that $\tilde n= 3k-1$ whence $m$
satisfies the conditions of Theorem \ref{main}. Thus, this series of
examples
shows that the value of $e_1$ can decrease
with respect to the general position
by any positive integer.

3. Let us generalize the previous example further.
Put \hfill\break $m=(m_1,m_2,1,\ldots,1)$ with an arbitrary $n$ satisfying
$m_1-m_2<n-2<m_1+m_2$.  (We still assume that
$\tilde n=m_1+m_2+n-2$ is even which is not very important).
It is easy to check that the conditions on $m_1,m_2$ and $n$ imply the
conditions (i) and (ii) of Theorem \ref{main} whence for an
arrangement in general position $e_1=\frac{\tilde n}{2}$.
The matrix $M$ for this example differs from a Vandermonde matrix
by the only one
gap in the degrees, down from $m_1-1$ to $\frac{m_1-m_2+n-4}{2}$
(after factoring out $z_i$ from the $i$th row for every $i$).
In other words $d_1=\pm s_{\lambda}$
where $s_{\lambda}$ is the Schur function for the rectangular
partition diagram with the base
$\frac{m_1+m_2-n}{2}$ and the height $\frac{-m_1+m_2+n-2}{2}$.
Thus we obtain the description of the set of arrangements with degenerate $e_1$
as the zero locus of the special  `rectangular' Schur function.

4. Let us recall Example 1 in this section. It has $\tilde n=8$
which is the smallest value of $\tilde n$ that can allow
different $e_1$. Besides the example 1, there is only one more
vector $m$ with $\tilde n=8$ such that $e_1$ is not determined by
$m$. It is $m=(3,2,2,1)$. Indeed $m$ satisfies the conditions (i)
and (ii) of Theorem \ref {main} and for an arrangement
$x^3y^2(x-z_1y)^2(x-z_2y)$  the condition for $e_1=3$ is
$d_1=-2z_1+z_2=0$.
This defines one (up to
projective isomorphism) arrangement with $z_1=1$ and $z_2=2$.
Notice that again the quadruple $(0,-2,\infty,-1)$ of points in
$\CP^1$ is harmonic (cf. the order of points from Example 1).

One more class of examples with a similar flavor is given by Proposition 7.4.
\bigskip
\section{The Terao conjecture}
\label{terao}

As we mentioned in the Introduction the main motivation for this work
is the Terao conjecture and Yoshinaga theorem for 3-arrangements.

\begin{conjecture} \textup{({\bf Terao, \cite{OT}})}
If two hyperplane arrangements over the same field
have isomorphic intersection lattices and
one of them is free then the other is free as well.
\end{conjecture}

This conjecture was posed (as a question) more than 20 years ago
(\cite{Te}) and is still open. The first non-trivial (and already hard)
case is formed by arrangements in $\CC^3$ or equivalently in $\CP^2$.
The main evidence in favor of the conjecture is a theorem by Terao
(\cite{OT}, Theorem 4.61) that
for a free arrangement its exponents are precisely the roots of
the characteristic polynomial of its intersection lattice
whence are determined by the intersection lattice.
We also mention a result from (\cite{Yu}) that in the space
of all arrangements with a fixed intersection lattice the free ones form a
 set open in the Zariski topology.

Some recent progress has been made by Yoshinaga in \cite{Yo}.
We paraphrase here the relevant part of his Theorem 2.2.

\begin{theorem} \textup{({\bf Yoshinaga, \cite{Yo}})}
If an arrangement $\A$ of lines in $\CP^2$ is free with exponents
${1,e_1,e_2}$ then its restriction to every $\ell\in\A$ has
(as a multi-arrangement)
exponents ${e_1,e_2}$. If the characteristic polynomial of $\A$
has roots $1,d_1,d_2$ and the restriction of $\A$ to some line has
exponents $d_1,d_2$ then $\A$ is free.
\end{theorem}

Yoshinaga's theorem implies in particular that if an arrangement $\A$ in
$\CP^2$ gives a counterexample to the Terao conjecture then its
restrictions to all lines allow simultaneous degeneration of their minimal
exponents. Since the multiplicity vectors of the restrictions are
determined
by the intersection lattice of $\A$, any time the multiplicity vector
on
one of the lines determines the exponents we obtain a class of
arrangements
(or lattices) satisfying the conjecture.

Let us give some examples to this effect.
To avoid any confusion
let us emphasize that the multiplicity $\bar m(p)$ of a point $p$ for an arrangement $\A$ in
$\CP^2$ is the number of all lines of the arrangement passing through $p$.
If $p\in\ell\in\A$ then the multiplicity of $p$
in the restricted to $\ell$ multi-arrangement is $m(p)=\bar m(p)-1$.

 Terao's conjecture holds for the following classes of arrangements $\A$ in $\CP^2$.

1. There exists a line in $\A$ whose intersections with other lines of
$\A$ are concentrated at no more than 3 points.

2. There exists a line in $\A$ such that all points on it
have multiplicities not larger than 3.

3. There exists a line in $\A$ such that the multiplicity
vector of the restricted multi-arrangement satisfies
$m_1\geq \sum\limits_{i=2}^nm_i$.

4. There exists a line in $\A$ such that the average of the
multiplicities of the restricted multi-arrangement
is less than 2.

\begin{remark}
\label{8}
If $\A$ is a line arrangement in $\CP^2$ such that $|\A|\leq 8$ then
either condition 1 or 4 holds. Hence, for $|\A|\leq 8$ Terao's conjecture
holds.
\end{remark}

In fact we can prove that a larger class of arrangements
satisfy Terao's conjecture.

\begin{proposition}
\label{fibertype}
Let $\A$ be a line arrangement in $\CP^2$. If there exists a intersection point $p$
such that $\bar m(p)>\frac{1}{2}(|\A|-3)$
then Terao's conjecture holds for $\A$ (i.e. for all arrangements with the intersection lattice isomorphic to that of $\A$).
\end{proposition}
\proof
Suppose $p$ satisfies the condition and denote by $\ell_1 ,\ell_2,...,\ell_k$ all the lines of $\A$ passing through $p$ where $k=\bar m(p)$. There are two alternatives.

(i) There are two lines $\ell,\ \ell'\in \A$ such that $p'=\ell \cap \ell'\notin \ell_i$ for every $i$. Then the points $\ell \cap \ell_1,...,\ell \cap \ell_k, p'$ are pairwise distinct. this implies that the number $k'$ of intersection points on $\ell$ satisfies $k'\geq k+1>\frac{1}{2}(|\A|-1)$ whence the average of the multiplicities of the restriction of $\A$ to $\ell$ is less than 2. The result follows from (4) above.

(ii) Every point of intersection lies on a line passing through $p$. This implies that the intersection lattice of $\A$ is supersolvable (see \cite{OT}, pp. 30,31). Since every arrangement with a supersolvable lattice is free (\cite{OT}, Theorem 4.58) the result follows. \qed

In particular, Remark \ref{8} generalizes to the following.

\begin{corollary}
\label{9,10}
Terao's conjecture holds for line arrangements in $\CP^2$ of cardinality less than 11. 
\end{corollary}
\proof
Indeed every arrangement with at most 10 lines either satisfies the condition 2 above or the condition of the previous proposition.                     \qed

\bigskip

\bibliographystyle{amsplain}

\begin{thebibliography}{99}

\bibitem{Mu}T.~Muir, A treatise on the theory of determinants, Revised and enlarged by William H. Metzler Dover Publications, Inc., New York 1960 vii+766 pp. 

\bibitem{OT} P.~Orlik and H.~Terao, Arrangements of hyperplanes, Springer
Verlag, 1992.

\bibitem{ST} L.~Solomon and H.~Terao, The double Coxeter arrangement, Commentarii Math. Helvetici 73 (1998), no. 2, 237-258.

\bibitem{Te} H.~Terao, The exponents of a free hypersurface,
Singularities,
editor P.~Orlik, Proc. Symp. Pure. Math 40, Pt. 2 (1983), 561-566.

\bibitem{Yo}  M.~Yoshinaga, Characterization of a free arrangement and conjecture of Edelman and Reiner, Invent. Math. 157 (2004), no. 2, 449-454.

\bibitem{Yu} S.~Yuzvinsky, Free and locally free arrangements with
a given intersection lattice, Proc. of AMS 118(1993), 745-752.

\bibitem{Zi} G.~Ziegler, Multi-arrangements of hyperplanes and their
freeness, in Singularities, Contemp. Math. 90(1986), 345-359.

\end{thebibliography}

\end{document}